\documentclass[11pt]{amsart} 
\usepackage{amsthm}
\usepackage[cyr]{aeguill }
\usepackage{dsfont }
\usepackage{fancyhdr} 
\usepackage{upgreek} 
\usepackage{mathrsfs} 
\usepackage{eufrak} 
\usepackage{stmaryrd}
\usepackage{amsmath, amsfonts, amssymb, color, subeqnarray} \usepackage{graphicx, longtable} 
\usepackage[T1]{fontenc} 
\usepackage[a4paper]{geometry} 
\usepackage[utf8]{inputenc}
\usepackage[english]{babel} 
\usepackage{cite}
\usepackage{datetime}
\theoremstyle{plain} 
\newtheorem{teo}{Theorem}[section]
\newtheorem{pro}{Proposition}[section]
 
\newtheorem{cor}{Corollary}[section] 
\theoremstyle{definition}
\newtheorem{defi}{Definition}[section] 
 
\theoremstyle{remark}
\newtheorem{nota}{\bf Remark}[section]

\def\sca #1_#2_#3 {\langle#1,#2\rangle_#3}

\def\C{\mathbb {C}}
\def\m{\mu}
\def\<{\langle}
\def\>{\rangle}

\def\n{\nu}

\def\pp{\varphi}
\def\H{{\mathcal H}}

\def\RR{{\Bbb X}}

\newcommand{\bq}{\begin{equation}}

\newcommand{\eq}{\end{equation}}
\newcommand{\ba}{\begin{eqnarray*}}
\newcommand{\ea}{\end{eqnarray*}}

\newcommand{\ove}[1]{\overline{#1}}

\newcommand{\ban}{\begin{eqnarray}}
\newcommand{\ean}{\end{eqnarray}}

\makeatletter
\g@addto@macro{\endabstract}{\@setabstract}
\newcommand{\authorfootnotes}{\renewcommand\thefootnote{\@fnsymbol\c@footnote}}%
\makeatother
\begin{document}
\date{\today}
\begin{center}
  {\large {\bf\uppercase{ On functions of negative type on the Olshanski spherical pair $(SL(\infty),SU(\infty))$}}} \par \bigskip

  \normalsize
	
  \authorfootnotes
  M. RABAOUI\let\thefootnote\relax\footnote{\today\par 
	\textsuperscript{b}Universit\'{e} de Tunis El Manar, Facult\'{e} des Sciences de Tunis, Laboratoire d'Analyse Math\'{e}matique et Applications, LR11ES11, 2092 El Manar Tunis, Tunisie.\par 
  \textsuperscript{c}Universit\'{e} de Carthage, Institut pr\'{e}paratoire aux \'{e}tudes d'ing\'{e}nieur de Nabeul, Campus Universitaire Merazka, 8000 Nabeul, Tunisie.\par }\textsuperscript{a,b,c} \par \bigskip

  \textsuperscript{a}{\small King Faisal University, College of Science, Department of Mathematics and Statistics,\\
	P.O. Box: 400, Al Ahsa 31982, Kingdom of Saudi Arabia }

	{\small  mrabaoui@kfu.edu.sa }

\end{center}

\bigskip

\begin{abstract} In this paper, using a generalized Bochner type representation for Olshanski spherical pairs, we prove the boundedness of every $SU(\infty)$-biinvariant continuous function of negative type on the infinite dimensional special linear group $SL(\infty)$.\vspace{0.3cm}
MS classification: \quad 22E30 (Primary); 43A35, 43A85, 43A90 (Secondary) 
\end{abstract}

\bigskip

\section{Introduction} 

\bigskip

The continuous functions of negative type have been introduced by Schoenberg (cf. \cite{Sho1}) in order to determine the invariant Hilbert metrics on a group. The importance of such functions lies in their applications in the theory of limit theorems for independent and identically distributed random variables (cf. \cite{Pet}). In 1990, G. Olshanski elaborated a spherical analysis theory related to the inductive limits of increasing sequences of Gelfand pairs (cf. \cite{Olsh}, \cite{Butov}, \cite{Far2} and \cite{Rab2}). In this new framework, many results concerning the continuous functions of negative type had been obtained (cf. \cite{BBouali}, \cite{Bouali} and \cite{Rab3}). \\ 


In this paper, we consider the spherical pair $(G_\infty, \ K_\infty)$, which is the inductive limit of the increasing sequence of Gelfand pairs $(G_n, \ K_n)$ where $G_n=SL(n,\mathbb C)$ and $K_n=SU(n)$ are respectively the special linear group and its unitary subgroup of order $n$. The infinite dimensional special linear group $G_\infty=SL(\infty)$ consists in unimodular infinite invertible complex matrices $(g_{ij})$ with a finite number of coefficients $g_{ij}\neq\delta_{ij}$ and $K_\infty=SU(\infty)$ is the group of unimodular infinite unitary matrices $(u_{ij})$ with complex coefficients such that $u_{ij}=\delta_{ij}$ for $i+j$ large enough. The group $G_\infty$ is equipped with the inductive limit topology and the subgroup $K_\infty$ is closed.\\


	
	
   A function $\pp:G_\infty\longrightarrow \mathbb C$ is said to be of positive type if the kernel defined on $G_\infty\times G_\infty$ by $(g_1,g_2)\longmapsto\pp(g_2^{-1}g_1)$ is of positive type, i.e. for all $g_1, g_2,\ldots, g_n \in G_\infty$ and all $c_1, c_2,\dots, c_n \in \C$, $$\sum_{i=1}^n\sum_{j=1}^nc_i\ove{c_j}\pp(g_j^{-1}g_i)\ \geq 0 .$$ Every function $\pp$ of positive type on $G_\infty$ is Hermitian, i.e. for all $g\in G_\infty$, $\ove{\pp(g)}=\pp(g^{-1})$. In addition, the function $\pp$ is bounded : $|\pp(g)|\leq\pp(e)$, where $e$ is the origin of $G_\infty$. Besides, a function $\pp$ defined on $G_\infty$ is said to be $K_\infty$-biinvariant if it holds that $\pp(k_1gk_2)=\pp(g)$, for all $k_1$, $k_2\in K_\infty$ and all $g\in G_\infty$. A function $\psi$, d{e}fined on $G_\infty$, with complex values is said to be {\it of negative type} if $\psi(e)\geq0,$ $\psi(g^{-1})=\ove{\psi(g)}$ and, for all $g_1,\dots,g_N\in G_\infty$ and all $c_1,\dots,c_N\in\Bbb C$ such that $\sum_{i=1}^Nc_i=0$, $$\sum_{i,j=1}^Nc_i\ove{c_j}\psi(g_j^{-1}g_i)\leq0.$$
	
    A spherical function for the Olshanski spherical pair $(G_\infty,K_\infty)$ is a $K_\infty$-biinvariant continuous function $\varphi$ on $G_\infty$ which is normalized at the origin of the group:  $\varphi(e) = 1$, and such that $$\lim_{n\rightarrow\infty}\int_{K_n}\varphi(xky)\alpha_n(dk)=\varphi(x)\varphi(y),$$ where $\alpha_n$  denotes the normalized Haar measure on $K_n$  (cf. \cite{Olsh} and \cite{Far8}). In section 2, using the explicit expression of the sph{e}rical functions of positive type relatively to the spherical pair $(G_\infty, \ K_\infty)$ given by N.I. Nessonov in \cite{Ness}, we prove that the spherical dual of the pair $(G_\infty, \ K_\infty)$ can be identified with a set $\RR$ of parameters. This enables us to obtain a parametrized version of the generalized Bochner theorem (\cite{Rab}, Theorem 7), which represents the key theorem used, in the last	section, to prove the boundedness of all continuous functions of negative type, relatively to the pair $(G_\infty, \ K_\infty)$. The method that we use follows an idea given by C. Berg, J.P. Christensen and P. Ressel in \cite{Ber1} and the work of M. Bouali in \cite{Bouali} and \cite{BBouali}. 

\section{Functions of class $ \mathfrak{B} $ : definition and convergence}


In his paper \cite{Ness}, N.I. Nessonov gives a complete classification of the representations of $ G_\infty $ containing the trivial representation of the unitary group $K_\infty$. This class of representations contains those of spherical representations. In fact, the trivial representation of $K_\infty$ is contained in a unitary representation $ (\pi, \H) $ of $ G_\infty $ if and only if the subspace $\H_{K_\infty}$ consisting of $K_\infty$-invariant vectors in $\H $ is not reduced to zero. On the other hand, every spherical representation is associated via the G.N.S. construction (Guelfand-Naimark-Segal) to a spherical function, i.e. an extreme element of the convex set of continuous functions of positive type on $G_\infty$ which are biinvariant by the subgroup $ K_\infty $ and equal to $1$ at the origin. Thus, this representation possesses a cyclic vector which is in addition $K_\infty$-invariant. It follows that the result of N.I. Nessonov, that one finds at the end of his paper (cf. \cite{Ness}, Theorem 4.5), allows us to give an explicit expression of spherical functions of positive type relatively to the pair $(G_\infty, K_\infty)$.

\bigskip

\begin{defi}\label{dif}\rm{
A {\it function of class $\mathfrak{B}$} of parameter $\;\alpha=(\alpha_1,\dots,\alpha_p)\in \mathbb R^p\;$ with $\;p\in\mathbb N\;$ is defined on $\;\mathbb R\;$ by :
$${\Pi}(\alpha,\lambda):=\prod_{j=1}^{p}\big(\cosh(\lambda)-i\alpha_j\sinh(\lambda)\big)^{-1}.$$}
\end{defi}

\bigskip

  A function $\Pi(\alpha,\lambda)$ of class $\mathfrak{B}$ is of positive type on $\mathbb R$. In fact, for $\alpha\in \mathbb R$, the function $\left(\cosh(\lambda)-i\alpha\sinh(\lambda)\right)^{-1}$ is the Fourier transform of the probability measure $\mu_\alpha$ on $\mathbb R$ given by $$\mu_\alpha=g_\alpha(t)\,dt \qquad g_\alpha(t)=\dfrac{e^{t\arctan \alpha}}{2\sqrt{1+\alpha^2}\cosh\big(\pi t/2\big)}.$$ This can be proved by the residue theorem via a contour integral of the function $$\left(\cosh(\pi z)-i\alpha\sinh(\pi z)\right)^{-1}e^{-2\pi itz} ,$$ around the rectangle of vertices $\pm R$, $\pm R+i$ and letting $R\to +\infty$. Then, the function of positive type $\Pi(\alpha,\lambda)$ is the Fourier transform of the convolution measure $\mu_{\alpha_1}\star\dots\star\mu_{\alpha_p}$, where $(\alpha_1,\dots,\alpha_p)\in \mathbb R^p$.

\bigskip

We consider, on the set of functions of class $ \mathfrak{B} $, the topology of uniform convergence on compact sets of $\mathbb R$. The resulting topological space is metrizable and complete. This topology can be expressed in terms of the set $\RR$ of parameters $\;\alpha=(\alpha_1,\dots,\alpha_p)\in \mathbb R^p\;$ with $\;p\in\mathbb N$. Let us introduce, for $\;p\in\mathbb N$ the sets $$\RR_p=\left\{(\alpha_1,\dots,\alpha_p)\in \mathbb R^p\; |\; \alpha_1\leq \alpha_2\leq\dots\leq\alpha_p\right\}.$$ Since, the class $\mathfrak{B}$ should contain the trivial function $1$, one has to define the space of one element $\RR_0=\{\varnothing\}$ with the discrete topology. In consequence, the space of parameters $\RR$ becomes the disjoint union of $\RR_p$ with $p=0,1,2,\dots$. For every $p\in \mathbb N_{_0}$, the set $\RR_p$ is both closed and open in $\RR$. Hence, it is locally compact with a countable basis.\\



 The set of functions of class $\mathfrak B$ is then parametrized by the set $\RR$ and every function $\Pi(\alpha,\lambda)$ have the following logarithmic derivative: 

\bq\label{logg}\begin{aligned}\frac{\Pi^{'}(\alpha,\lambda)}{\Pi(\alpha,\lambda)}&=\, i\sum\limits_{m=0}^{\infty}\widetilde{p}_m(\alpha)\left(i\tanh(\lambda)\right)^{m},\end{aligned}\eq  where $p_m$ is the Newton power sum function : for $\alpha=(\alpha_1,\alpha_2,\ldots,\alpha_p)\in\mathbb R^p$ and $m\in\mathbb N_0$, $$p_m(\alpha)=\sum_{k=1}^{p}\alpha_k^m.$$

 \bq\label{powernewton}\widetilde{p}_m(\alpha)=p_{m+1}(\alpha)+p_{m-1}(\alpha) \quad (m\in\mathbb N_0),\eq  with the convention that  $$p_{-1}(\alpha)=0 \ .$$

\bigskip

\begin{pro}\label{bij} The mapping $\alpha\longmapsto {\Pi}(\alpha,.)$ is a bijection of $\;\RR^{^*}=\RR\setminus\RR_0\;$ to the set of functions of class $\mathfrak B$.\\\end{pro} 
\begin{Proof} For $\lambda\in\mathbb R$, $\alpha\in\RR_p$ and $\beta\in\RR_q$, let us assume that

 $$\Pi(\alpha,\lambda)=\Pi(\beta,\lambda). $$ Since,

 $$\lim_{\lambda\to+\infty}e^{p\lambda}\, \Pi(\alpha,\lambda)=2^p\prod_{k=1}^p(1-i\alpha_k)^{-1},\vspace{0.3cm}$$ it follows that $$p=q.\vspace{0.3cm}$$ Using the power series expansion (\ref{logg}), we conclude that $\,\alpha=\beta$.$\hspace{1cm}\Box\\$ 
\end{Proof}

A function of class $ \mathfrak{B} $ of parameter $\alpha$ is the Fourier transform of a probability measure $\m_\alpha$. Let $\mathfrak M_\RR$ be the set of these measures : $$\mathfrak M_\RR=\left\{\m_\alpha  | \Pi(\alpha,\lambda)=\widehat{\m_\alpha}(\lambda)\right\}.$$

We consider on $\mathfrak M_\RR$ the weak topology of measures. We will prove that the topology defined on $\RR$ is equivalent to the weak topology of $\mathfrak M_\RR$.\\ 


\begin{pro}\label{irman} The topology of $\RR$ is equivalent to the topology of $\mathfrak M_\RR$.
\end{pro}

\begin{Proof} (i) Assume that $\alpha^{(n)}$ converges to $\alpha$ in the topology of $\RR$. Then $\alpha\in\RR_p$, for some $p$ and so, there exists $n_0$ such that $\alpha^{(n)}\in\RR_p$, for $n\geq n_0$. It follows that $\Pi(\alpha^{(n)},\lambda)$ converges uniformly on compact sets of $\mathbb R$ to $\Pi(\alpha,\lambda)$. Finally, by applying the L{\'e}vy-Cramer theorem, one can prove that $\m_{\alpha^{(n)}}$ converges weakly to $\m_\alpha$.\\

(ii) Assume that $\m_{\alpha^{(n)}}$ converges weakly to $\m_\alpha$. This implies that $\Pi(\alpha^{(n)},\lambda)$ converges uniformly on compact sets in $\mathbb R$ to $\Pi(\alpha,\lambda)$. Let $\lambda_0>0$. Since, the function $\Pi(\alpha^{(n)},\lambda)$ is continuous, non-zero on $\mathbb R$ and satisfies $\Pi(\alpha^{(n)},0)=1$, there exists $0<C<1$ and $n_1\in\mathbb N$ such that, for every $n\geq n_1$,

 $$\Re(\Pi(\alpha^{(n)},\lambda_0))\geq C,$$
and so,

$$|\Pi(\alpha^{(n)},\lambda_0)|\geq C.\vspace{0.3cm}$$ 
In consequence,

$$\prod_{j=1}^{p_n}\big(\cosh^2(\lambda_0)+{\alpha_j^{(n)}}^2\sinh^2(\lambda_0)\big)^{-1}\geq 
C^{2}.\vspace{0.3cm}$$ Hence, we get 

$$\sinh^2(\lambda_0)\sum_{j=1}^{p_n}{\alpha_j^{(n)}}^2\leq 
\prod_{j=1}^{p_n}\big(1+{\alpha_j^{(n)}}^2\sinh^2(\lambda_0)\big)\leq\frac{1}{C^{2}},$$ or $$\label{1}\sum_{j=1}^{p_n}{\alpha_j^{(n)}}^2\leq \frac{1}{C^{2}\sinh^2(\lambda_0)}.\vspace{0.3cm}$$ This implies that, for all $n\geq n_1$ and all $1\leq j\leq p_n$, 

\bq\label{bound1}|\alpha_j^{(n)}|\leq \frac{1}{C\sinh(\lambda_0)}.\vspace{0.3cm}\eq Moreover, 

 $$\prod_{j=1}^{p_n}\cosh^2(\lambda_0)\leq\frac{1}{C^{2}}.\vspace{0.3cm}$$ As a result, for all $n\geq n_1$,

 \bq \label{bound2}p_n\leq-\frac{\log C}{\log(\cosh(\lambda_0))}.\vspace{0.3cm}\eq 

Since, by (\ref{bound2}), the sequence $p_n$ is bounded, we can assume, upon passing to a subsequence, that $p_n = p$. Further, by (\ref{bound1}), there exists a subsequence $n_k$ such that $\alpha^{(n_k)}$ converges to some $\beta\in\RR_p$. It follows from the first part of the proof that the corresponding sequence of measures $\mu_{\alpha^{(n_k)}}$ converges weakly to $\mu_\beta$. Thus, we get $\mu_\alpha=\mu_{_\beta}$ and so $\alpha=\beta$. Since any accumulation point of the sequence $\alpha^{(n)}$ in $\RR_p$ is necessary equal to $\alpha$, it follows that the sequence $\alpha^{(n)}$ itself converges  to $\alpha$ in $\RR$. $\hspace{1cm}\Box$\\

\end{Proof}


\section{Parametrized Bochner Theorem on $(G_\infty, \ K_\infty)$}

\bigskip

Let $\mathscr{P}$ be the convex set of $K_\infty$ biinvariant continuous functions of positive type on $G_\infty$. The topology defined on the set of its extremal points ${\rm ext(}\mathscr{P}{\rm )}$ can be expressed in terms of the set $\RR$. This enables us to prove a parametrized version of the generalized Bochner theorem (\cite{Rab}, Theorem 7).\\

\begin{cor}\label{tti} The mapping $\alpha\longmapsto {\Pi}(\alpha,.)$ is a homeomorphism between $\RR^{^{*}}$ and $\mathfrak B$.\\
\end{cor}
\begin{Proof} By Proposition \ref{bij}, the function $\alpha\mapsto \Pi(\alpha,.)$ is a bijection between $\RR^{^{*}}$ and $\mathfrak B$. Besides, by Proposition \ref{irman}, the topology of uniform convergence on compact sets of $\mathbb R$ defined on $\mathfrak{B}$ is equivalent to the topology defined on the set of parameters $\RR^{^{*}}$. In consequence, the application $\alpha\mapsto {\Pi}(\alpha,.)$ defines a homeomorphism between $\RR$ and $\mathfrak B$.$\hspace{1cm}\Box$\\
\end{Proof}

\bigskip

The group $G_\infty$ can be seen as the group of unimodular infinite invertible matrices $g=[g_{ij}]_{i,j=1}^\infty$ with a finite number of coefficients $g_{ij}\neq\delta_{ij}$. Every matrix $g$ in $G_\infty$ can be written under the form $g=u$ diag($e^{\lambda_1}, e^{\lambda_2},\ldots$) $v$, with $u,v\in K_\infty$, $\lambda_1,
\lambda_2,\ldots\in\Bbb R$, $\sum_j\lambda_j=0$ and $\lambda_j = 0$ for $j$ large enough. We say that a $K_\infty$-biinvariant function $\pp$ on $G_\infty$ which is normalized by the condition $\pp(e)=1$, is multiplicative if
$$\pp(\rm{diag}\it (a_{\rm{1}}, \it a_{\rm{2}}, \ldots))=\rm\Phi(\it a_{\rm{1}})\rm\Phi\it(a_{\rm{2}})\ldots,$$ where $\Phi$ is a function defined on $\Bbb R$ such that $\Phi(0)=1$. In other words, for all $g\in
G_\infty$, $$\pp(g)=\Phi(a_1)\Phi(a_2)\ldots,$$ where $a_1,a_1,\ldots\in\Bbb R$ are the elements of the diagonal matrix in the decomposition of $g$.\\
\begin{teo}{\rm(\cite{Far2}, Theorem 5.2, page 228)} Let $\pp$ be a $K_\infty$-biinvariant continuous function of positive type on $G_\infty$ such that $\pp(e)=1$. Then, $\pp$ is spherical (or extremal) if and only if it is multiplicative.\\
\end{teo}

 By the preceding theorem, the set ext($\mathscr P$) consists of the multiplicative functions of $\mathscr P$. The following theorem gives these elements via the functions of class $\mathfrak{B}$.\\

\begin{teo}\label{bi}
The extremal points of the set $\mathscr P$ are the functions $\varphi_\alpha$ defined for all $g=u$ \rm{diag(}\it$e^{\lambda_1}, e^{\lambda_2},\ldots$\rm{)}\it $v\in G_\infty$, with $u,v\in K_\infty$, $\lambda_1,
\lambda_2,\ldots\in\Bbb R$, $\sum_j \lambda_j=0$ and $\lambda_j = 0$, for $j$ large enough, by
$$\pp_\alpha(g)=\pp_\alpha\big({\rm diag(}e^{\lambda_1}, e^{\lambda_2},\ldots{\rm)}\big)=\prod_{j=1}^\infty\Pi(\alpha,\lambda_j).$$

\begin{paragraph}{\bfseries{Proof.}}
{\rm ${\rm\bf a)}$ By the preceding theorem, a function in $\mathscr P$ defined, for all $g\in
G_\infty$, by $\varphi_\alpha (g)=\prod_j\Pi(\alpha,\lambda_j)$ is extremal in $\mathscr P$ since it is multiplicative.\\
${\rm\bf b)}$ Let $\varphi\in {\rm ext}(\mathscr P)$. It is a spherical function. By the result of Nessonov, the function $\pp$ is a product of functions of class $\mathfrak{B}$. $\hspace{1cm}\Box$\\}
\end{paragraph}
\end{teo}

\begin{nota} {\rm  Remark that the trivial spherical function is isolated in ${\rm ext}(\mathscr P)$. This follows immediately from the fact that there is no sequence $\alpha^{(n)}$ for which $\Pi(\alpha^{(n)},.)$ converges to $1$. It is a well-known fact (cf. \cite{Far5}, Theorem 5.2) that, in the case of a Gelfand pair $(G,K)$, the isolation of the trivial spherical function in the spherical dual implies the boundedness of every $K$-biinvariant continuous function of negative type on $G$. This fact is not true for the inductive limit groups since, in general, they are not locally compact .

 }
\end{nota}

\begin{pro} The correspondence $\RR\leftrightarrow$ {\rm ext(}$\mathscr P${\rm)} is an isomorphism between two standard spaces.\\
\begin{paragraph}{\bfseries{Proof.}}{\rm Since the set $\RR$ is locally compact, separable, metrizable and complete, it represents, in consequence, a standard space. In addition, the proof of the generalized Bochner theorem (\cite{Rab}, Theorem 7), shows that ext($\mathscr P$) is a Borel subset of a standard space. Hence, it is standard by (\cite{D2}, Appendix B, B 20). Furthermore, the correspondence ext($\mathscr P$)$\rightarrow\RR$, $\pp_\alpha\mapsto \alpha$ is Borelian and one-to-one. In consequence, by (\cite{D2}, Appendix B, B22), it is an isomorphism between two standard spaces.$\hspace{1cm}\Box$}
\end{paragraph}
\end{pro}

Hence, we can get a parametrized version of the generalized Bochner theorem (\cite{Rab}, Theorem 7):

\bigskip

\begin{teo}\label{th3} Let $\varphi$ be a $K_\infty$-biinvariant continuous function of positive type on $G_\infty$. Then, there exists a unique positive and bounded measure $\mu$ defined on $\RR$ such that, for every $g\in G_\infty$,
$$\varphi(g)=\int_{\RR}\varphi_{\alpha}(g)\mu(d\alpha).$$

\end{teo}


\bigskip

 The following proposition follows from the uniqueness of the representing measure in the generalized Bochner Theorem.\\

\begin{pro}\label{uniq} Let $\mu_1$ and $\mu_2$ be two positive and bounded measures on $\RR$ satisfying  \begin{equation}\label{on}\int_{\RR}(1-\pp_\alpha(g))\,\mu_1(d\alpha)=\int_{\RR}(1-\pp_\alpha(g))\,\mu_2(d\alpha).\end{equation} Then,

$$\mu_1=\mu_2\quad\mbox{on}\quad\RR^{^*}. $$


\end{pro}

\bigskip

\begin{paragraph}{\bfseries{Proof.}}\rm{Since the function $\pp_\alpha$ is spherical, it satisfies, for all $g,h\in G_\infty$,

\begin{equation}\label{0}\lim_{n\to\infty}\int_{K_n}\pp_\alpha(gkh)\,dk=\pp_\alpha(g)\pp_\alpha(h),\end{equation}
where $dk$ is the normalized Haar measure of the group $K_n$. Since, the measures $dk$ and $\mu_1(d\alpha)$ are bounded, the function $$\displaystyle(\alpha, k)\mapsto |1-\pp_\alpha(gkh)|$$
is int{e}grable with respect to the product measure $\displaystyle\mu_1(d\alpha)\times dk$. Hence, by the Fubini theorem, we get that

 $$\int_{K_n}\int_{\RR}(1-\pp_\alpha(gkh))\,\mu_1(d\alpha)\,dk=\int_{\RR}\int_{K_n}(1-\pp_\alpha(gkh))\,dk\,\mu_1(d\alpha).\vspace{0.5cm}$$ Using the fact that $|1-\pp_\alpha(gkh)|\leq 2$ and that the measure $dk$ is a probability measure, we conclude that the function

$$\alpha\mapsto\int_{K_n}(1-\pp_\alpha(gkh))\,dk$$ is dominated by $2$. Since, the measure $\displaystyle\mu_1(d\alpha)$ is positive and bounded, the dominated convergence theorem and the equation (\ref{0}), imply that

$$\lim_{n\to\infty}\int_{\RR}\int_{K_n}(1-\pp_\alpha(gkh))\,dk\,\mu_1(d\alpha)=\int_{\RR}(1-\pp_\alpha(g)\pp_\alpha(h))\,\mu_1(d\alpha),$$ and so,

$$\lim_{n\to\infty}\int_{K_n}\int_{\RR}(1-\pp_\alpha(gkh))\,\mu_1(d\alpha)\,dk=\int_{\RR}(1-\pp_\alpha(g)\pp_\alpha(h))\,\mu_1(d\alpha).\vspace{0.5cm}$$ In consequence, by equation (\ref{on}), we get that

$$\int_{\RR}(1-\pp_\alpha(g)\pp_\alpha(h))\,\mu_1(d\alpha)=\int_{\RR}(1-\pp_\alpha(g)\pp_\alpha(h))\,\mu_2(d\alpha).\vspace{0.5cm}$$ By substituting, in the previous equation, $g$ by $g^{-1}$, and using the fact that $\pp_\alpha$ is Hermitian, we get

$$\int_{\RR}(1-\pp_\alpha(h)\Re\pp_\alpha(g))\,\mu_1(d\alpha)=\int_{\RR}(1-\pp_\alpha(h)\Re\pp_\alpha(g))\,\mu_2(d\alpha).\vspace{0.5cm}$$ By substituting $g$ by $h$ in equation (\ref{on}), and considering the difference with the previous equation, we get, for all $g,h\in G_\infty$,

$$\int_{\RR}\pp_\alpha(h)(1-\Re\pp_\alpha(g))\,\mu_1(d\alpha)=\int_{\RR}\pp_\alpha(h)(1-\Re\pp_\alpha(g))\,\mu_2(d\alpha).\vspace{0.5cm}$$  Now, let us consider the function $\widetilde\varphi$ defined on $G_\infty$ by

$$\widetilde\varphi(h)=\int_{\RR}\pp_\alpha(h)(1-\Re\pp_\alpha(g))\,\mu_1(d\alpha).\vspace{0.5cm}$$ Since, for all $g\in G_\infty$, we have $1-\Re\pp_\alpha(g)\geq 0$, the function $\widetilde\varphi$ is of positive type. It is also $K_\infty$-biinvariant and continuous on $G_\infty$. Moreover, the function $\widetilde\varphi$ is bounded, since it is dominated by a $\mu_1$-integrable function on $\RR$ : 

$$|\pp_\alpha(h)(1-\Re\pp_\alpha(g))|\leq 1-\Re\pp_\alpha(g).\vspace{0.5cm}$$ So, for all $g\in G_\infty$, $$\int_{\RR}\pp_\alpha(h)\,\mu_{1,g}(d\alpha)=\int_{\RR}\pp_\alpha(h)\,\mu_{2,g}(d\alpha).\vspace{0.5cm}$$
where,
$$\mu_{i,g}(d\alpha)=(1-\Re\pp_\alpha(g))\,\mu_i(d\alpha)
\quad \rm{for}\quad \it i=\rm 1, 2.\vspace{0.5cm}$$ By uniqueness of the representing measure in the generalized Bochner theorem (Theorem \ref{th3}), we get that, for all $g\in G_\infty$,

$$\mu_{1,g}=\mu_{2,g}\quad\mbox{on}\quad\RR.\vspace{0.5cm}$$ Since, we have $|\Pi(\alpha,\lambda)|<1$, for $\alpha\in\RR^{^*}$ and $\lambda\neq0$, we get, for $g_0={\rm diag}(e^1,e^{-1},1,...)$, that $1-\Re\pp_{\alpha}(g_0)>0$. Hence

$$\mu_1=\mu_2\quad\mbox{on}\quad\RR^{^*}. \qquad\qquad\hfill\square$$}


\end{paragraph}

\section{Functions of negative type on $(G_\infty,K_\infty)$}

\bigskip

 In this section, using the generalized Bochner theorem (Theorem \ref{th3}), we establish the boundedness of every continuous function of negative type on the pair $(G_\infty,K_\infty)$. The method that we follow is inspired from the work of M. Bouali in \cite{Bouali} and \cite{BBouali}. If $\varphi$ is a function of positive type, then $\psi(g)=\varphi(e)-\varphi(g)$ is a bounded function of negative type. The functions of negative type and those of positive type are related by the following property :\\

\begin{pro}\label{pr}{\rm(Schoenberg (\cite{Sho1}, page 527) and (\cite{Ber2}, Theorem 7.8))}\\ The function $\psi$ is of negative type if and only if 
$\psi(e)\geq 0,$ and, for all $t\geq 0$, $e^{-t\psi}$ is of positive type.
\end{pro}

\bigskip


 


\begin{teo}\label{th2} Every $K_\infty$-biinvariant continuous function $\psi$ of negative type on 
$G_\infty$ is bounded and can be uniquely written as $\displaystyle \psi(g)=\psi(e)+\varphi(e)-\varphi(g)$, where $\varphi$ is a function of positive type relatively to the spherical pair $(G_\infty,K_\infty)$.\\

\begin{paragraph}{\bfseries{Proof.}}\rm{ Let $\psi$ be a $K_\infty$-biinvariant continuous function of negative type on $G_\infty$. Since $\psi(g)-\psi(e)$ is also $K_\infty$-biinvariant, continuous and of negative type, we can assume, without loss of generality, that $\psi(e)=0$. For $t\geq 0$, the function $e^{-t\psi}$ is $K_\infty$-biinvariant, continuous and of positive type on $G_\infty$. Hence, by the generalized Bochner theorem (Theorem \ref{th3}), there exists a probability measure $\mu_t$ on $\RR$ such that

 $$e^{-t\psi(g)}=\int_{\RR}\varphi_{\alpha}(g)\;\mu_t(d\alpha).\vspace{0.5cm}$$ It follows that, for all $t> 0$,

\begin{equation}\label{rw}e^{-t\Re\psi(g)}\cos(t\Im\psi(g))=\int_{\RR}\Re\pp_\alpha(g)\
  \m_t(d\alpha),\end{equation}
and
$$e^{-t\Re\psi(g)}\sin(t\Im\psi(g))=-\int_{\RR}\Im\pp_\alpha(g)\m_t(d\alpha).\vspace{0.3cm}$$

${\rm\bf i)}$ The equation (\ref{rw}) implies that
 
$$\frac{1-e^{-t\Re\psi(g)}\cos(t\Im\psi(g))}{t}=\int_\RR(1-\Re\pp_\alpha(g))\frac{\m_t}{t}(d\alpha).$$
In addition, we have $$\lim_{t\to
  0}\frac{1-e^{-t\Re\psi(g)}\cos(t\Im\psi(g))}{t}=\Re\psi(g),$$  
$$\lim_{t\to
  +\infty}\frac{1-e^{-t\Re\psi(g)}\cos(t\Im\psi(g))}{t}=0.\vspace{0.3cm}$$ For $g$ fixed, the last expression is a continuous function in $t$ on $]0,+\infty[$ which tends to $0$ as $t$ tends to $\infty$. Thus, there exists a constant $C(g)\geq 0$ such that
	
$$0\leq\frac{1-e^{-t\Re\psi(g)}\cos(t\Im\psi(g))}{t}\leq C(g).$$
Therefore,
$$\int_\RR(1-\Re\pp_\alpha(g))\frac{\m_t}{t}(d\alpha)\leq
C(g).\vspace{0.3cm}$$
In particular, for $g_0={\rm diag}(e^1,e^{-1},1,...)$, we have

\begin{equation}\label{eq1}\int_\RR(1-\Re\pp_\alpha(g_0))\frac{\m_t}{t}(d\alpha)\leq 
C(g_0)=M.\vspace{0.3cm}\end{equation} Let $\n_t$ denote the positive and bounded measure defined on $\RR$ by

$$\n_t=(1-\Re\pp_\alpha(g_0))\frac{\m_t}{t}.\vspace{0.3cm}$$
Since the set $\left\{\n_t\;|\;t>0\right\}$ is uniformly bounded by the constant $M$, it is relatively compact for the weak topology $\sigma({\mathcal M}(\RR),\mathscr
C_0(\RR))$, where ${\mathcal M}(\RR)$ is the set of positive and bounded measures on $\RR$ and $\mathscr C_0(\RR)$ is the set of continuous functions on $\RR$, vanishing at $\infty$. In consequence, there exists a sequence $(t_j)$ in $]0, +\infty[$ converging to $0$, such that the measures $\n_{t_j}$ weakly converge to a positive and bounded measure $\n$, i.e. for all $f\in {\mathscr C}_0(\RR)$,

 $$\lim_{j\rightarrow\infty}\int_{\RR}f(\alpha)\,\n_{t_j}(d\alpha)=
\int_{\RR}f(\alpha)\,\n(d\alpha).\vspace{0.3cm}$$ Since $\pp_\alpha(g_0)=1$ for $\alpha_0=\{\varnothing\}$, the measure $\n_t$ has no mass at the point $\alpha_0$ and therefore $\n\left(\{\varnothing\}\right)=0$. On another hand, we have 

\bq\label{eqo}\begin{split}
\frac{1-e^{-t_j\Re\psi(g)}\cos(t_j\Im\psi(g))}{t_j}&=\int_{\RR^{^*}}\left[\frac{1-\Re\pp_\alpha(g)}{1-\Re\pp_\alpha(g_0)}-1\right]\n_{t_j}(d\alpha)\\&+\frac{1-e^{-t_j\Re\psi(g_0)}\cos(t_j\Im\psi(g_0))}{t_j}.\end{split}\eq

For $g\neq e$ fixed, let us consider the function $f$ defined on $\RR^{^*}$ by

$$f(\alpha)=\frac{1-\Re\pp_\alpha(g)}{1-\Re\pp_\alpha(g_0)}-1.\vspace{0.3cm}$$ 
 The function $f$ is well defined since $\Re\pp_\alpha(g_0)\neq 1$, for all $\alpha\in\RR^{^*}$. It is also continuous, by Corollary \ref{tti} and Theorem \ref{bi}. In addition, for every $g\in G_\infty$ such that $\lambda_j\neq0$ for at least one $j$, we have
$$0<|\pp_\alpha(g)|\leq\prod_{j=1}^\infty\prod_{k=1}^p\big(\cosh^2(\lambda_j)+\alpha_k^2\sinh^2(\lambda_j)\big)^{-\frac{1}{2}}.\vspace{0.3cm}$$ It follows that, for $\alpha\in\RR_p$ and $p_0$ large enough independent of $\alpha$,  $$|\pp_\alpha(g)|\leq \cosh^{-p}(\lambda_j)<\varepsilon.\vspace{0.3cm} $$ Besides, for $1\leq p \leq p_0$ and for $\alpha$ in $\RR_p$ such that $||\alpha||$ is large enough,  $$|\pp_\alpha(g)|\leq \dfrac{1}{||\alpha||\,|\sinh(\lambda_j)|}<\varepsilon.\vspace{0.3cm} $$ In consequence, $|\pp_\alpha(g)|$ tends to $0$ as $\alpha$ tends to $\infty$ in the locally compact space $\RR$. Therefore, the function $f$ belongs to ${\mathscr C}_0(\RR^{^*})$. Now, as $j$ tends to $\infty$ in the equation (\ref{eqo}), one gets, for all $g\neq e$,

$$\begin{aligned}\Re\psi(g)=\int_{\RR^{^*}}f(\alpha)\,\n(d\alpha)+\Re\psi(g_0)=\int_{\RR^{^*}}\left(\frac{1-\Re\pp_\alpha(g)}{1-\Re\pp_\alpha(g_0)}-1\right)\n(d\alpha)+\Re\psi(g_0).\end{aligned}\vspace{0.3cm}$$
 Since $$\left|\frac{1-\Re\pp_\alpha(g)}{1-\Re\pp_\alpha(g_0)}-1\right| \leq \frac{2}{1-\Re\pp_\alpha(g_0)}+1,\vspace{0.3cm} $$ with $\Re\pp_\alpha(g_0)\leq\cosh^{-2p}(1) $, the integrand in the last integral is dominated and we can apply the dominated convergence theorem to get, as $g$ tends to $e$, 

$$\displaystyle\n(\RR^{^*})=\Re\psi(g_0).\vspace{0.3cm}$$ So, $$\Re\psi(g)=\int_{\RR^{^*}}(1-\Re\pp_\alpha(g))\mu(d\alpha),\vspace{0.3cm}$$
where $\mu$ is the measure defined on $\RR^{^*}$ by 
$$\mu=\frac{1}{1-\Re\pp_\alpha(g_0)}\n.$$ By Corollary \ref{tti}, the function
$$\alpha\mapsto\frac{1}{1-\Re\pp_\alpha(g_0)}$$ is continuous on $\RR^{^*}$. Moreover, this function tends to $1$ as $\alpha$ tends to $\infty$ and it has as limit $\cosh^{-2p}(1)$ at $0$. It follows that it is bounded and so is the measure $\mu$.\\

${\rm\bf ii)}$ Similarly, we prove that the function $h$ defined on $\RR^{^*}$ by
$$h(\alpha)=\frac{\Im\pp_\alpha(g)}{1-\Re\pp_\alpha(g_0)}\ ,$$
belongs to ${\mathscr C}_0(\RR^{^*})$. As $j$ tends to $\infty$ in the relation\\
$$e^{-t_j\Re\psi(g)}\frac{\sin(t_j\Im\psi(g))}{t_j}=-\int_{\RR^{^*}}\frac{\Im\pp_\alpha(g)}{1-\Re\pp_\alpha(g_0)}\n_{t_j}(d\alpha),\vspace{0.5cm}$$
we get $$\begin{aligned}\Im\psi(g)=-\int_{\RR^{^*}}\frac{\Im\pp_\alpha(g)}{1-\Re\pp_\alpha(g_0)}\n(d\alpha)=-\int_{\RR^{^*}}\Im\pp_\alpha(g)\mu(d\alpha).\end{aligned}\vspace{0.5cm}$$ Finally, 
$$\psi(g)=\int_{\RR^{^*}}(1-\pp_\alpha(g))\mu(d\alpha) = \pp(e)-\pp(g),\vspace{0.15cm}$$ where $$\pp(g)= \int_{\RR^{^*}}\varphi_{\alpha}(g)\,\mu(d\alpha).$$   Remark that, by Proposition \ref{uniq}, the measure $\mu$ is unique.   \qquad \qquad $\hfill\square$}
\end{paragraph}
\end{teo}



\begin{nota} {\rm  In \cite{Kaz}, Kazhdan introduced the notion of the property (T) as follows: a locally compact topological group $G$ has the property (T) if the trivial representation ${\bf{\iota}}_{G}$ is isolated in the  unitary dual $\widehat{G}$ which is the set of all unitary irreducible representations of $G$ equipped by the Fell topology. The notion of property (T) was extended to any Hausdorff topological group and was also related to other topological properties such as the property (FH). One can refers to \cite{Bek} for extensive reading. Property (T) is related to the functions of negative type since it holds that if a Hausdorff topological group $G$ has property (T), then every continuous function of negative type on $G$ is bounded. The converse, which holds in the locally compact case, is not true in general. A question immediately arises : does $SL(\infty)$ have property (T)?\\
 }
\end{nota}

\begin{paragraph}{\bfseries{Acknowledgment}}{\rm The author would like to thank the referee for his valuable comments that improved considerably the exposition of the results in this paper. The author is also grateful to the Deanship of Scientific Research at King Faisal University for financially supporting this work under Project 160273.\\}

\end{paragraph}

\bibliographystyle{amsalpha}

\end{document}